\documentclass{amsart}
\usepackage{amsmath}
\usepackage{amssymb}
\usepackage[dvips]{graphicx}
\usepackage[latin1]{inputenc}
\usepackage{amssymb,latexsym}

\newtheorem{thm}{Theora}[section]
\newtheorem{Theorem}[thm]{Theorem}

\newtheorem{Lemma}[thm]{Lemma}
\newtheorem{Proposition}[thm]{Proposition}
\newtheorem{Definition}[thm]{Definition}
\newtheorem{Definitions}[thm]{Definitions}
\newtheorem{Remark}[thm]{Remark}

\newtheorem{Definition and notations}[thm]{Definition and notations}
\newtheorem{Proof}[thm]{Proof}
\title[Existence of weak solutions for nonlinear ...]
{Existence of weak solutions for nonlinear elliptic systems involving the $(p(x),q(x))$-Laplacian}

\author[ Mounir Hsini ]{Mounir Hsini }

\begin{document}
\maketitle

\medskip
{\footnotesize \centerline{Institut Préparatoire aux Etudes d'Ingénieurs de Tunis}
\smallskip
}

\medskip
{\footnotesize\centerline{mounir.hsini@ipeit.rnu.tn}
\medskip
\bigskip
}
\begin{abstract}
In this paper, we prove the existence of weak solutions for the following nonlinear elliptic
system
\begin{eqnarray*}
\begin{array}{lll}
-\Delta_{p(x)}u \; = \; a(x)|u|^{p(x)-2}u \; - \; b(x)|u|^{\alpha(x)}|v|^{\beta(x)} v \; + \; f(x)\;\;\; in\;\;\Omega, \\
-\Delta_{q(x)}v \; = \; c(x) |v|^{q(x)-2}v \; - \;d(x)|v|^{\beta(x)}|u|^{\alpha(x)}  u \; + \; g(x)\;\;\; in \;\;\Omega,\\
 u \; = \; v \; = \; 0 \quad on \;\;\partial\Omega,
\end{array}
\end{eqnarray*}
where $\Omega$ is an open bounded domains of $\mathbb{R}^N$ with a smooth boundary $\partial\Omega$
and $\Delta_{p(x)}$ denotes the $p(x)$-Laplacian.The existence of weak solutions is proved using
the theory of monotone operators. Similar result will be proved when $\Omega=\mathbb{R}^N$.
\end{abstract}
\section{Introduction}
The purpose of this paper is to study the existence of weak solutions for the following nonlinear elliptic system involving the $p(x)$-Laplacian.
\begin{eqnarray}\label{S1}
\left\{\begin{array}{lll}
-\Delta_{p(x)}u \; = \; a(x)|u|^{p(x)-2}u \; - \; b(x)|u|^{\alpha(x)}|v|^{\beta(x)} v \; + \; f(x)\;\;\; in\;\;\Omega,\\
-\Delta_{q(x)}v \; = \; c(x) |v|^{q(x)-2}v \; - \;d(x) |v|^{\beta(x)}|u|^{\alpha(x)} u \; + \; g(x)\;\;\; in \;\;\Omega,\\
 u \; = \; v \; = \; 0 \quad on \;\;\partial\Omega,
\end{array}\right.
\end{eqnarray}
where $\Omega$ is an open bounded domains in $\mathbb{R}^N$ with a smooth boundary
$\partial\Omega$.  The operator $\Delta_{p(x)}u = \mbox{div}\big{(} \mid \nabla u\mid^{p(x) - 2} \nabla u  \big{)}$
is called  $p(x)$-Laplacian, which will be reduced to the $p$-Laplacian when $p(x)=p$ a constant.\\
The study of various mathematical problems with variable exponent has been received considerable attention in recent years, for examples we cite works of X-L Fan and V. Radulescu \cite{F}, \cite{MR}. The operator $p(x)$-Laplacian turns up in many mathematical settings, e.g., Non-Newtonian fluids, reaction-diffusion problems, porous media, astronomy, quasi-conformal mappings..etc. see  \cite{a3,a4,d1}.\\ Problems including this operator for bounded domains have been studied in \cite{F, MR} and for unbounded domains in \cite{EADE, HB, f4a}.
Many authors have studied semilinear and non linear elliptic systems, as a reference we cite \cite{CST, EADE, MH, HV, PV}.\\

The generalized formulation for many stationary boundary value problems for partial differential equations leads to operator equation of type
$$L(u)=f$$
on a Banach space. Indeed, the weak formulation consists in looking
for an unknown function $u$ from a Banach space $H$ such that an
integral identity containing $u$ holds for each test function $v$
from the space $H$. Since the identity is linear in $v$, we can take
its sides as values of continuous linear functionals at the element
$v\in H$. Denoting the terms containing unknown $u$ as the value of
an operator $A$, we obtain
$$
(L (u),v)=(f,v)\quad \forall v\in H,
$$
which is equivalent to equality of functionals on $H$, i.e. the
equality of elements of $H'$ (the dual space of $H$): $L (u)=f$.\\
In this paper, we consider nonlinear systems with model $L$ of the
form
\begin{eqnarray*}
 L\{u,v\}&=&\{-\Delta_{p(x)}u-a(x)|u|^{p(x)-2}u+b(x)|u|^{\alpha(x)}|v|^{\beta(x)} v\;,\\
 & &\qquad\qquad - \;\Delta_{q(x)}v + c(x)|u|^{\alpha(x)} |v|^{\beta(x)} u - d(x) |v|^{q(x)-2}v \}.
\end{eqnarray*}
When $p(x)=p$ is constant, the existence of solutions for such systems was proved,
using the method of sub and super solutions in \cite{b3,b4,s1}.
In this study, we use another technique for proving the  existence of weak
solutions. We need the theory of monotone operators.\\
To resolve the system $(\ref{S1})$, we introduce the following intermediary problem
\begin{eqnarray}\label{P4}
\left\{\begin{array}{lll}
-\Delta_{p(x)} u= a(x)|u|^{p(x)-2}u+f(x), \quad x \in \Omega,\\
 u=0\qquad\text{on } \partial\Omega,
\end{array}\right.
\end{eqnarray}
where $\Omega$ is a bounded domain of $\mathbb{R}^N$, $p(.)\in
C^0(\overline{\Omega})$ satisfying $\underset{x\in \Omega}{\inf}\,p(x)>1$
and $\gamma$ is a non negative function in $L^\infty(\Omega).$
\begin{Theorem}\label{T1}
The nonlinear elliptic problem $(\ref{P4})$ have a non trivial weak solution.
\end{Theorem}
\begin{Theorem}\label{T2}
Under the (so-defined) assumptions F0), F1), F2) and F3).
The nonlinear elliptic system (\ref{S1}) have a non trivial weak solution, when $\Omega$ is a bounded domain or $\Omega=\mathbb{R}^N$.
\end{Theorem}
This paper consists of five sections. First, we recall some elementary proprieties
of the Generalized Lebesgue-Sobolev Spaces and introduce the notations needed in this work. Section 3 is devoted to the study of some preliminary results
which allows us to prove the existence of weak solutions of our problem. Particulary we give the proof of the first Theorem. In the fourth section, we justify the existence of weak solutions in the case of bounded domains. The goal of the last section is the main result, when $\Omega=\mathbb{R}^N$.
\section{Generalized Lebesgue-Sobolev Spaces Setting.}
In order to discuss problem (\ref{S1}), we need some theories on
spaces  $W^{1,p(x)}(\Omega)$ which we call generalized Lebesgue-
Sobolev spaces.  Let us shortly recall some basic facts about the
setup for generalized Lebesgue- Sobolev spaces, for more details see
for instance \cite{F}, \cite{Julian Musielak}, \cite{DZF} and
\cite{DZQ}.\\
Let
\begin{eqnarray*}
C_+(\overline{\Omega})= \{ h\;/\; h\in C(\overline{\Omega}), h(x) > 1 \text{for any }x\in \overline{\Omega} \}.
\end{eqnarray*}
For $p(x)\in C_+(\overline{\Omega})$, we define the
variable exponent Lebesgue space $ L^{p(x)}(\Omega)$ by
$$ L^{p(x)}(\Omega)=\{u\;/\;\; u \;\text{is a measurable real-valued
function},\; \int_\Omega|u(x)|^{p(x)} dx < \infty\}.$$ We define the
so-called Luxemburg norm, on this space by the formula
\begin{eqnarray*}
|u|_{L^{p(x)}}=\inf\{\alpha>0,\;\,\int_\Omega\Big{|}
\frac{u(x)}{\alpha}\Big{|}^{p(x)}dx\leq1\}.
\end{eqnarray*}
It's well known, that $(L^{p(x)}(\Omega); |.|_{L^{p(x)}})$ is a is a separable, uniform
convex Banach space. \\$(L^{p(x)}(\Omega); |.|_{L^{p(x)}})$  is
termed a generalized Lebesgue space. Moreover, its conjugate space
is $L^{p'(x)}(\Omega)$, where $\frac{1}{p'(x)}+\frac{1}{p(x)}=1$. For
$u\in L^{p(x)}(\Omega)$ and $v\in L^{p'(x)}(\Omega)$, one has the
following inequality
\begin{eqnarray}\label{3}
\Big{|}\int_\Omega u(x)v(x) dx\Big{|}\leq
\Big{(}\frac{ 1}{p^-}+\frac{1}{p'^-}\Big{)}|u|_{L^{p(x)}}|v|_{L^{p'(x)}}\leq 2|u|_{L^{p(x)}}|v|_{L^{p'(x)}},
\end{eqnarray}
where, $p^-=\underset{\overline{\Omega}}{\min}\, p(x)$ and
$q^-=\underset{\overline{\Omega}}{\min}\, q(x)$.\\ Note that
$L^{p_{2}(x)}(\Omega)\hookrightarrow L^{p_{1}(x)}(\Omega)$, for every
functions $p_{1}$ and $p_{2}$ in $C(\overline{\Omega})$ satisfying
$p_{1}(x)\leq p_{2}(x),$ for any $x\in \overline{\Omega}$. In addition
this imbedding is continuous. \\
An important role in manipulating the generalized Lebesgue spaces is
played by the modular of the $L^{p(x)}(\Omega)$ space, which is the
mapping $\rho_{p(x)}:\,L^{p(x)}(\Omega)\rightarrow \mathbb{R}$
defined by
$$ \rho_{p(x)}(u)=\int_\Omega |u|^{p(x)} dx.$$
If $(u_n), u \in L^{p(x)}(\Omega)$ and $p^+<\infty$ then the
following relations hold true.
\begin{eqnarray}\label{e3}
|u|_{L^{p(x)}}>1 \Rightarrow |u|_{L^{p(x)}}^{p^-}\leq
\rho_{p(x)}(u) \leq |u|_{L^{p(x)}}^{p^+},
\end{eqnarray}
\begin{eqnarray}\label{e4}
|u|_{L^{p(x)}}< 1 \Rightarrow |u|_{L^{p(x)}}^{p^+}\leq
\rho_{p(x)}(u) \leq |u|_{L^{p(x)}}^{p^-},
 \end{eqnarray}
\begin{eqnarray}\label{e5}
|u_n-u|_{L^{p(x)}}\rightarrow0\;\;\mbox{if and
if}\;\,\rho_{p(x)}(u_n-u)\rightarrow0.
\end{eqnarray}
Another property interesting the variable exponent Lebesgue space $ L^{p(x)}(\Omega)$ is
\begin{Proposition} \label{p22} $($ see \cite{Ed}$)$\;
Let $p(x)$ and $q(x)$ be measurable functions such that $p(x)\in L^\infty(\mathbb{R}^N)$ and $1\leq p(x)q(x)\leq \infty,$ for a.e. $x\in \mathbb{R}^N$ Let $u\in  L^{q(x)}(\mathbb{R}^N)$, $u\neq 0$. Then
$$|u|_{p(x)q(x)}\leq 1\Rightarrow |u|_{p(x)q(x)}^{p^+}\leq ||u|^{p(x)}|_{q(x)}\leq |u|_{p(x)q(x)}^{p^-},$$
$$|u|_{p(x)q(x)}\geq 1\Rightarrow |u|_{p(x)q(x)}^{p^-}\leq ||u|^{p(x)}|_{q(x)}\leq |u|_{p(x)q(x)}^{p^+}.$$
In particular, if $p(x)=p$ is constant, then
$$||u|^{p}|_{q(x)}=|u|_{pq}^p.$$
\end{Proposition}
The generalized Lebesgue-Sobolev space is defined by:
$$W^{1,p(x)}(\Omega)=\{u\in L^{p(x)}(\Omega)\;\;\mbox{such that}\;
\;|\nabla u|\in L^{p(x)}(\Omega)\}.$$ $W^{1,p(x)}(\Omega)$ can be
equipped with the norm defined as follow
\begin{eqnarray}\label{6}
\|u\|_{p(x)}= |u|_{L^{p(x)}}+|\nabla
u|_{L^{p(x)}},\;\mbox{for all}\;\;u\in W^{1,p(x)}(\Omega).
\end{eqnarray}
In this paper, we denote by $W_0^{1,p(x)}(\Omega)$ the closure of
$C_0^\infty(\Omega)$ in $W^{1,p(x)}(\Omega)$.\\
Due to Fan and Zhao \cite{F}, generalized Lebesgue- Sobolev spaces
$W^{1,p(x)}(\Omega)$ and $W_0^{1,p(x)}(\Omega)$ are separable
reflexive Banach spaces. On the other hand if $q\in
C_+(\overline{\Omega})$ satisfying $q(x)<p^*(x)$ for any $x\in
\overline{\Omega}$, the imbedding from $W^{1,p(x)}(\Omega)$ into
$L^{q(x)}(\Omega)$ is compact and continuous. Note that Poincar\'e
inequality is also satisfied and we have existence of a constant
$C>0$ such that
\begin{eqnarray}\label{e6}
|u|_{L^{p(x)}}\leq C  |\nabla u|_{L^{p(x)}}, \;\,\mbox{for
all}\;u\in W_0^{1,p(x)}(\Omega).
\end{eqnarray}
In view of (\ref{6}), it follows that $|\nabla u|_{L^{p(x)}}$ and
$\|u\|_{p(x)}$ are equivalent norms on $W_0^{1,p(x)}(\Omega)$.
Hence, we will use $|\nabla u|_{L^{p(x)}}$ to replace
$\|u\|_{p(x)}$.
\begin{Definitions}\label{d0}
$1<p(x)<N$ and for $x\in \mathbb{R}^N$, let define
$$
p^*(x)=
\left\{\begin{array}{ll}
\frac{Np(x)}{N-p(x)}\;\;\,p(x)<N,\\
+\infty \;\;\;\;\;\,\;\,p(x)>N,
\end{array}\right.
$$
where $p^*(x)$ is the so-called critical Sobolev exponent of $p(x)$.
\end{Definitions}
\begin{Proposition} \label{p222} $($ see \cite{Ed}$)\;$
Let $p(x)\in C^{0,1}_+(\mathbb{R}^N),$ that is Lipshitz-continuous function defined on
$\mathbb{R}^N$, then there exists a positive constant $c$ such that $$|u|_{p^*(x)}\leq \|u\|_{p(x)},$$
for all $u\in W_0^{1,p(x)}(\Omega)$.
\end{Proposition}
Through this paper we suppose that the following assumptions are satisfied.\\
$$\mbox{F0)}\hspace{4.5cm} a(x), c(x)\;\; \mbox{are resp. in}\;\; L^{p'(x)}(\Omega)\mbox{and}\;\;L^{q'(x)}(\Omega).$$
$$\mbox{F1)}\hspace{1cm}s(x)=\frac{p(x)p^*(x)q^*(x)}{p(x)p^*(x)q^*(x)-pq^*(x)-p^*(x)q^*(x)},\hspace{1cm} b(x)\in L^{s(x)}(\Omega),$$
$$\mbox{F2)}\hspace{1cm} r(x)=\frac{q(x)p^*(x)q^*(x)}{q(x)p^*(x)q^*(x)-qq^*(x)-p^*(x)q^*(x)}, \hspace{1cm} d(x)\in L^{r(x)}(\Omega),$$
$$\widetilde{p}(x)=\frac{p(x)p^*(x)}{p^*(x)-p(x)};\hspace{1cm}
\widetilde{q}(x)=\frac{q(x)q^*(x)}{q^*(x)-q(x)},$$
$$\mbox{F3)}\hspace{1cm}1<p^-, q^-\hspace{0.2cm} \alpha^+<p^--1, q^--1, \hspace{0.2cm}\beta^+<p^--1, q^--1,p^+<p^*>2.$$
Others notations will be introduced as we need.
\section{Preliminary}
This section is devoted to the study of problems of type: $Lu=f$, where $L$ is an operator from $H$ (Banach space) into it's dual $H^*$. The tools needed for such aim is the variational method, more precisely theory of monotone operator.\\
To this end, we introduce some technical results \cite{b2,b4,z1} which allows us to the proof of Theorem\ref{T1}. Note that hypothesis F0) and F3) will be used in this section.
First, we recall the following definition.
\begin{Definitions}\label{d1}
Let $L : H \to H^*$  be an operator on a Banach space $H$.\\ We say that the
operator $L $ is:\\
\textmd{1.}  Monotone  if $\langle L (u_{1}) - L (u_{2}),u_{1}-u_{2}\rangle \ge 0$ for all $u_{1}, u_{2}$.\\
\textmd{2.}  Strongly continuous  if $u_n \rightharpoonup u$  implies $L (u_n) \to L (u)$.\\
\textmd{3.}   Weakly continuous  if $u_n \rightharpoonup u$  implies $L (u_n) \rightharpoonup L (u)$.\\
\textmd{4.}  Demi-continuous   if $u_n \to u$  implies $L (u_n) \rightharpoonup L (u)$.\\
\textmd{5.} The operator $L $ is said to satisfy the   $M_0$-condition  if $ u_n\rightharpoonup u$,  $L (u_n) \rightharpoonup f$ and $\langle L (u_{n}),u_{n}\rangle \to \langle f,u \rangle$   imply $ L (u)=f$.
\end{Definitions}
The following Proposition plays an important role in the present paper. Precisely, it gives a sufficient conditions to the existence of weak solutions for the problems $Lu=f$.
\begin{Proposition}\label{p0}
Let $H$ be a separable reflexive Banach space and $L : H\to H'$
an operator which is: coercive, bounded, demicontinuous,
and satisfying $M_0$ condition.
Then the equation $L (u)=f$ admits a solution for each $f\in H'$.
\end{Proposition}
Next, we consider the eigenvalue problem involving the $p(x)$-Laplacian of the form
\begin{eqnarray}\label{p3}
\left\{\begin{array}{lll}
-\Delta_{p(x)} u=\lambda\; a(x)|u|^{p(x)-2}u, \quad x \in \Omega,\\
 u=0\qquad\text{on } \partial\Omega,
\end{array}\right.
\end{eqnarray}
where $\Omega$ is a bounded domain of $\mathbb{R}^N$, $p(.)\in
C^0(\overline{\Omega})$ satisfying $\underset{x\in \Omega}{\inf}\,p(x)>1$
and $\,a(x)$ is a non negative function in $L^\infty(\Omega).$\\
Below we write $X=W_0^{1,p(.)}$ and $\|u\|=|\nabla u|_{p(x)}.$
\begin{Definition}
Let $\lambda\in \mathbb{R}$ and $u\in X.$ $(u,\lambda)$ is called a solution of problem ($\ref{p3}$) if
$$\int_\Omega|\nabla u|^{p(x)-2}\nabla u\nabla v dx= \lambda\int_\Omega a(x)|u|^{p(x)-2}uv dx;\;\,\,\forall\, v\in X.$$
\end{Definition}
If $(u,\lambda)$ is a solution of problem ($\ref{p3}$) and $u\in X\backslash\{0\}$, we say $\lambda$ and $u$
an eigenvalue and an eigenfunction corresponding to $\lambda$ for problem ($\ref{p3}$), respectively. We recall that an eigenvalue $\lambda$ is called principal if there exists a nonnegative eigenfunction corresponding to $\lambda$, i.e., if there exists a nonnegative $u\in X\backslash\{0\}$ such that $(u,\lambda)$ is a solution of ($\ref{p3}$). \\
Now, we are ready to introduce a technical Lemma which is a consequence of Theorem 3.8 of X. Fan, for a reference, we cite \cite{F2}.
\begin{Lemma}\label{thm4}
Under assumptions above, Problem ($\ref{p3}$) has a solution ($\lambda_1(\,a(x)); u)$ satisfying $\lambda_1(\,a(x))>p^+/p^-$.
\end{Lemma}
In this sequel, we introduce the operator $L$ defined on $W^{1,p(x)}(\Omega)$ by
\begin{eqnarray}\label{Lu}
L u= -\Delta_{p(x)} u - a(x)|u|^{p(x)-2}u.
\end{eqnarray}
where $a(x)$ is a non negative function in $L^\infty(\Omega)$ and $\Omega$ is a bounded domain
of $\mathbb{R}^N$.\\
In order to prove the existence of weak solutions of the problem \ref{Lu}, we will need variational method. Precisely, we justify that the operator $L$ satisfies hypothesis of Proposition \ref{p0}. To this end, we introduce a series of Lemmas dealing with continuity, boundness, coercivity and monotonicity.  First we deal with continuity and boundness.
\begin{Lemma}\label{thm6}
$L $ is a bounded and demicontinuous operator.
\end{Lemma}
{\bf Proof. } It's clear that $L$ is the sum of $L_1$ and $L_2$, where
$$(L_1(u),v)=\int_\Omega|\nabla u|^{p(x)-2}\nabla u \nabla v dx\;\;\mbox{and }\;\,(L_2(u),v)=\int_\Omega a(x)|u|^{p(x)-2}u v dx.$$
Let us first prove the demicontinuity of the operator $L_1$.\\
Let $(u_n)\subset W^{1,p(x)}_{0}(\Omega)$ such that $u_n\rightarrow u_n$ in $W^{1,p(x)}_{0}(\Omega)$. We pass to a subsequence and assume that $u_n\rightarrow u$ and $\nabla u_n\rightarrow \nabla u$ pointwise almost everywhere. By the continuity of the map $\xi\mapsto |\xi|^{p(x)-2}\xi$, it follows that $|\nabla u_n|^{p(x)-2}\nabla u_n\rightarrow|\nabla u|^{p(x)-2}\nabla u$ almost everywhere. Since
$$\int_\Omega|\nabla u_n|^{p(x)-2}\nabla u_n|^{p(x)/(p(x)-1)}\,dx =\int_\Omega|\nabla u_n|^{p(x)}\,dx\leq M<\infty.$$
by the convergence of the sequence $(u_n)$, $(|\nabla u_n|^{p(x)-2}\nabla u_n)$ is bounded in $L^{p'(x)}(\Omega)$. Thus we may pass to a further sequence and assume that $(|\nabla u_n|^{p(x)-2}\nabla u_n)\rightarrow |\nabla u|^{p(x)-2}\nabla u$ weakly in $L^{p'(x)}(\Omega)$. This implies that the whole sequence converges weakly. Indeed: assuming the opposite, we find a weak neighbourhood $U$ of $|\nabla u|^{p(x)-2}\nabla u$ and a subsequence  such that $(|\nabla u_{n_k}|^{p(x)-2}\nabla u_{n_k})$ in $W^{1,p'(x)}_{0}(\Omega)$. We may assume pointwise convergence by passing to a further subsequence, and this sub-subsequence converges weakly in $L^{p'(x)}(\Omega)$ to $|\nabla u|^{p(x)-2}\nabla u$ by the earlier argument, which is a contradiction. It follows that
$$(L_1(u_n),v)=\int_\Omega|\nabla u_n|^{p(x)-2}\nabla u_v \nabla v dx\rightarrow \int_\Omega|\nabla u|^{p(x)-2}\nabla u \nabla v dx=(L_1(u_n),v).$$
therefore the demicontinuity of $L_1$.\\
Denotes
$$\Omega_1=\{x\in \Omega;\,|u(x)|\geq 1\}\;\;\mbox{and}\;\;\;\Omega_2=\{x\in \Omega;\,|u(x)|<1\},$$
Then
$$(L_2(u),v)=\int_{\Omega_1} a(x)|u|^{p(x)-2}u v dx+\int_{\Omega_1} a(x)|u|^{p(x)-2}u v dx.$$
In view of assumption $p^+-1<p^*(x)$, the following embeddings hold true:
$$W^{1,p(x)}_{0}(\Omega)\hookrightarrow L^{(p^+-1)p(x)}(\Omega) \;\;\mbox{and}\;\;
W^{1,p(x)}_{0}(\Omega)\hookrightarrow L^{(p^--1)p(x)}(\Omega).$$
Due Proposition \ref{p22}, we obtain
\begin{equation}\label{e20}
||u|^{p^+-1}|_{p(x)}=|u|^{p^+-1}_{kp(x)}\leq c_1\|u\|^{p^+-1}_{p(x)}.
\end{equation}
Take the function $a(x)$ in $L^{p^*(x)/(p^*(x)-2)}(\Omega),\;|u|^{p^+-1}, v\in L^{p^*(x)}(\Omega),$ and applying Holder inequality, we get
$$\Big{|}\int_{\Omega_1} a(x)|u|^{p(x)-2}u v dx\Big{|}\leq c_1 |\,a(x)|_{p^*/(p^*-2)}||u|^{p^+-1}|_{p^*(x)}|v|_{p^*(x)}$$
$$\hspace*{3.8cm}\leq c_2|\,a(x)|_{p^*/(p^*-2)}|u|^{p^+-1}_{(p^+-1)p^*(x)}|v|_{p^*(x)}.$$
$$\hspace*{4.0cm}\leq c_3|\,a(x)|_{p^*/(p^*-2)}\|u\|^{p^+-1}_{p(x)}\|v\|_{p(x)}<\infty.$$
Similarly,
$$\Big{|}\int_{\Omega_2} a(x)|u|^{p(x)-2}u v dx\Big{|}\leq c_4 |\,a(x)|_{p^*/(p^*-2)}\|u\|^{p^--1}_{p(x)}\|v\|_{p(x)}<\infty.$$
It follows that the operator $L_2(u,v)$ is well defined and bounded. Consequently $L$ is a bounded operator.
The proof of the demicontinuity of $L_2$ will be deduced from the following assumptions.\\
\textit{First step.} For all $u, v\in L^{p(x)}(\Omega), \; |u-v|_{p(x)}\rightarrow 0\Rightarrow ||u-v|^{\frac{p(x)}{p(x)-1}}|_{p(x)-1}\rightarrow 0$.\\
Let $\varepsilon>0$, $\eta < \varepsilon$ and
 $u, v\in L^{p(x)}(\Omega)$ such that $|u-v|_{p(x)}<\eta,$ then we have $$|u-v|_{p(x)}=\inf\{\mu\in ]0,\eta[;\;\int_\Omega\frac{|u-v|^{p(x)}}{\mu^{p(x)}}\;dx \leq 1\}.$$
On the other hand $<\mu<\eta<1$, it follows
$$\int_\Omega\frac{|u-v|^{p(x)}}{\mu^{p(x)-1}}\;dx\leq \int_\Omega\frac{|u-v|^{p(x)}}{\mu^{p(x)}}\;dx$$
and consequently
$$\inf\{\mu\in ]0,\eta[;\;\int_\Omega\frac{|u-v|^{p(x)}}{\mu^{p(x)-1}}\;dx \leq 1\}\leq \inf\{\mu\in ]0,\eta[;\;\int_\Omega\frac{|u-v|^{p(x)}}{\mu^{p(x)}}\;dx \leq 1\}.$$
Since the last term of this inequality represent $|u-v|_{p(x)}<\eta<\epsilon.$ The proof of the first claim will be immediately deduced if we consider the fact
$$||u-v|^{\frac{p(x)}{p(x)-1}}|_{p(x)-1}\leq \inf\{\mu\in ]0,\eta[;\;\int_\Omega\frac{|u-v|^{p(x)}}{\mu^{p(x)-1}}\;dx \leq 1\}.$$
\textit{Second step.} We claim that the map $u\in L^{p(x)}(\Omega)\rightarrow |u|^{p(x)-2}\in L^{\frac{p(x)}{p(x)-1}}(\Omega)$ is continuous. To this end we will use the convention
\begin{eqnarray*}
u^{p(x)}=
\left\{\begin{array}{ll}
 u^{p(x)}, \;  \;\;\mbox{for}\;\;\;u\geq 0; \\
-(-u)^{p(x)}\; \; \;\mbox{for}\;\; \;u\leq 0.
\end{array}\right.
\end{eqnarray*}
Our intention is to show the following identity:
$$|u-v|_{p(x)}\rightarrow0 \Rightarrow |u^{p(x)-1}-v^{p(x)-1}|_{\frac{p(x)}{{p(x)-1}}}\rightarrow0.$$
The result is trivial when $p(x)=2.$ We claim to prove the result for $p(x)>2$.
$$\rho_{\frac{p(x)}{{p(x)-1}}}(u^{p(x)-1}-v^{p(x)-1}):=\int_\Omega
|u^{p(x)-1}-v^{p(x)-1}|^{\frac{p(x)}{{p(x)-1}}}\;dx,$$
then, for $x\in \Omega$, by Lagrange theorem applied to the function $g(y)=y^{p(x)-1}$, there exists $c(x)$ somewhere between $u(x)$ and $v(x)$ satisfying
$$\frac{g(u(x))-g(v(x))}{u(x)-v(x)}=g'(c(x)).$$
Due to the fact that $|u-v|\in L^{p(x)}(\Omega)$, we have $|u-v|^{\frac{p(x)}{{p(x)-1}}}\in L^{p(x)-1}(\Omega)=(L^{\frac{p(x)-1}{p(x)-2}}(\Omega))^*$ and $|u|, |v|\in L^{p(x)}(\Omega)$ imply $|u|^{\frac{p(x)(p(x)-2)}{{p(x)-1}}},\;|v|^{\frac{p(x)(p(x)-2)}{{p(x)-1}}} \in L^{\frac{p(x)-1}{p(x)-2}}(\Omega).$ Hence
$$\rho_{\frac{p(x)}{{p(x)-1}}}(u^{p(x)-1}-v^{p(x)-1})\leq p^{+\frac{p^+}{p^--1}}
\int_\Omega |u-v|^{\frac{p(x)}{{p(x)-1}}}\sup(|u|, |v|)^{\frac{p(x)(p(x)-2)}{{p(x)-1}}}\,dx.$$
Thus the proof of the continuity by using  (\ref{3}), (\ref{e5}) and the second claim.\\
Ours second tools Lemma deals with coercivity, precisely we have
\begin{Lemma}\label{thm5}
The operator $L $ is coercive.
\end{Lemma}
{\bf Proof. } Let $\lambda_1(\,a(x))$ the first eigen value of the problem
$$-\Delta_{p(x)}u=\lambda a(x) |u|^{p(x)-2}u.$$
It's useful to recall the variational characterization
$$\lambda_1(\,a(x))=\inf\big{\{} \frac{\int_{\Omega}1/p(x)|\nabla u|^{p(x)}dx}{\int_{\Omega}a(x)/p(x)|u|^{p(x)}dx};\;\, u\in W_0^{1,p(x)}\backslash\{0\}\big{\}}.$$
Hence
$$\lambda_1 (\,a(x))\int_{\Omega}a(x)/p(x)|u|^{p(x)}dx\leq \int_{\Omega}1/p(x)|\nabla u|^{p(x)}dx,$$
and
$$\frac{\lambda_1(\,a(x))}{p^+}\int_{\Omega}a(x)|u|^{p(x)}dx\leq 1/p^-\int_{\Omega}|\nabla u|^{p(x)}dx.$$
It yields
\begin{equation}\label{e70}
\int_{\Omega}a(x)|u|^{p(x)}dx\leq \frac{p^+}{\lambda_1(\,a(x))p^-}\int_{\Omega}|\nabla u|^{p(x)}dx.
\end{equation}
On the other hand the operator $L$ satisfies
\begin{equation}\label{e99}
( L u , u)=\int_\Omega|\nabla u|^{p(x)}dx-
\int_\Omega a(x)| u|^{p(x)}dx.
\end{equation}
Combining equations (\ref{Lu}), (\ref{e70}) and (\ref{e99}), we obtain
$$ ( L u , u)\geq  \Big(1-\frac{p^+}{\lambda_{1}(\,a(x))p^-}\Big)\int_\Omega|\nabla u|^{p(x)}dx=(1-\frac{p^+}{\lambda_1(\,a(x))p^-})\rho_{p(x)}(|\nabla u|).$$
In view of Proposition \ref{p22} and Lemma \ref{thm4}, we obtain
$$( L u , u)\geq \inf(|\nabla u|_{p(x)}^{p^-}, |\nabla u|_{p(x)}^{p^+})
= \inf(\|u\|_{p(x)}^{p^-}, \|u\|_{p(x)}^{p^+}).$$
Using the fact that $p^->1$, one writes
$$( L u, u)/{\|u\|}\geq \inf(\|u\|_{p(x)}^{p^--1}, \|u\|_{p(x)}^{p^+-1})\to
\infty\quad\text{as}\quad \|u\|\to \infty.$$
Hence, the operator $L$ is coercive as required.\\
The third technical result in this section deals with monotonicity, in particular
\begin{Lemma}\label{thm7}
The operator $L $ is strictly monotone.
\end{Lemma}
{\bf Proof.} For the convenience, we give the idea of the proof. Recall the following elementary inequalities \cite{Ki;V} and \cite{Th}, from which we can get the strictly monotonicity of the operator $L$.
\begin{equation}\label{e13}
2^{2-p}|a-b|^p\leq \Big{(}a|a|^{p-2}-b|b|^{p-2}\Big{)}.(a-b),\;\,\mbox{if}\;\, p(x)\geq 2,
\end{equation}
\begin{equation}\label{e14}
(p-1)|a-b|^2\Big{(}|a|+|b|\Big{)}^{p-2}\leq \Big{(}a|a|^{p-2}-b|b|^{p-2}\Big{)}.(a-b),\;\,\mbox{if}\;\, 1<p(x)< 2.
\end{equation}
for all $a,b\in \mathbb{R}^n$, where .
denotes the standard inner product in $\mathbb{R}^n$.\\
\textbf{Remark.} Using previous Lemmas, all conditions of Proposition \ref{p0} are fulfilled. hence, the proof of Theorem 1 is completed.
\section{Nonlinear systems on bounded domains}
The goal of this section is to prove  existence of weak solutions for the system
\begin{eqnarray*} \label{P}
(\mathcal{S})\left\{\begin{array}{lll}
-\Delta_{p(x)}u \; = \; a(x)|u|^{p(x)-2}u \; - \; b(x)|u|^{\alpha(x)}|v|^{\beta(x)} v \;
+ \; f(x)\;\;\; in\;\;\Omega ,\\
-\Delta_{q(x)}v \; = \; c(x) |v|^{q(x)-2}v \; - \;d(x)|v|^{\beta(x)}|u|^{\alpha(x)}  u \;
+ \; g(x)\;\;\; in \;\;\Omega ,\\
 u \; = \; v \; = \; 0 \quad on \;\;\partial\Omega ,
\end{array}\right.
\end{eqnarray*}
where $\Omega$ is a bounded domain of $\mathbb{R}^{n}$, $p(x)$ and $q(x)$ are Lipshitz-continuous functions defined on $\mathbb{R}^N$. In addition, we suppose that $p(x), q(x)\in C^{0,1}(\Omega)$. We denote by $p'(x), q'(x)$ the conjugate exponent of $p(x), q(x)$ respectively. i.e.
$$\frac{1}{p(x)}+\frac{1}{p'(x)}= \frac{1}{q(x)}+\frac{1}{q'(x)}=1.$$
$ a(x), b(x), c(x), d(x)$ are non negative
functions satisfying condition F0), F1) and F2). Finally,  $\alpha(x)$ and $\beta(x)$ are regular nonnegative functions such that the assumption F3) will be satisfied.\\
In the following discussions, we will use the product space
\begin{equation}\label{e140}
W_{p(x),q(x)}:=W_0^{1,p(x)}(\Omega)\times W_0^{1,q(x)}(\Omega),
\end{equation}
which is equipped with the norm
\begin{equation}\label{e141}
\|(u,v)\|_{p(x),q(x)}:=\max\{\|u\|_{p(x)};\|v\|_{q(x)}\};\;\forall\,(u,v)\in W_{p(x),q(x)},
\end{equation}
where $\|u\|_{p(x)}$ (resp., $\|u\|_{q(x)}$) is the norm of $W_0^{1,p(x)}(\Omega)$ (resp., $W_0^{1,q(x)}(\Omega)$).\\
The space $W^*_{p(x),q(x)}$ denotes the dual space of $W_{p(x),q(x)}$ and equipped with the norm
$$\|.\|_{*,p(x),q(x)}:=\|.\|_{*p(x)}+\|.\|_{*,q(x)},$$
where $\|.\|_{*p(x)}, \|.\|_{*,q(x)}$ are respectively the norm of $W_0^{-1,p'(x)}(\Omega)$ and $W_0^{-1,q'(x)}(\Omega)$, dual resp.
of $W_0^{1,p(x)}(\Omega)$ and $W_0^{1,p(x)}(\Omega)$.\\
At beginning, we recall the following definition.
\begin{Definition}
$(u,v)\in W_{p(x),q(x)}$ is called a weak solution of the system $(\mathcal{S})$, if
$$\int_{\Omega}(|\nabla u|^{p(x)-2}\nabla u \nabla\Phi_1 +|\nabla v|^{p(x)-2}\nabla v\nabla\Phi_2 )dx = \int_{\Omega}(F_1(x,u,v)\Phi_1 +F_2(x,u,v)
\Phi_2 )dx,$$
for all $(\Phi_{1},\Phi_{2})\in W_{p(x),q(x)},$ where $F$ and $G$ are defined by
$$F_1(x,u,v)=a(x)|u|^{p(x)-2}u - b(x)|u|^{\alpha(x)}|v|^{\beta(x)} v +f(x),$$
$$F_2(x,u,v)=c(x)|v|^{q(x)-2}u - d(x)|u|^{\alpha(x)}|v|^{\beta(x)}u +g(x).$$
\end{Definition}
\begin{Remark}\label{R2}
The weak formulation of the system $(\mathcal{S})$ is reduced to the
operator form identity
\begin{eqnarray}
L_1 (u,v)- L_2 (u,v)+B(u,v)=F,
\end{eqnarray}
where $L_1, L_2, B$ and $F$ are defined on $W_{p(x),q(x)}$ as follow:
\begin{eqnarray*}
( L_1 (u,v), (\Phi_{1},\Phi_{2})) &=& \int_\Omega|\nabla u|^{p(x)-2}\nabla
u \nabla \Phi_{1}dx+ \int_\Omega|\nabla v|^{q(x)-2}\nabla v \nabla
\Phi_{2}dx,\\
( L_2 (u,v), (\Phi_{1},\Phi_{2}))&=&\int_\Omega a(x)
|u|^{p(x)-2}u\Phi_{1}dx+\int_\Omega  c(x) |v|^{q(x)-2}v \Phi_{2}dx,\\
(B(u,v), (\Phi_{1},\Phi_{2}))&=&
 \int_\Omega b(x) |u|^{\alpha(x)} |v|^{\beta(x)} v \Phi_{1}dx+
\int_\Omega d (x) |v|^{\beta(x)} |u|^{\alpha(x)} u \Phi_{2}dx,\\
(F,\Phi) &:=& ((f, g),(\Phi_{1},\Phi_{2}))=\int_\Omega  f
\Phi_{1}\;dx +\int_\Omega  g \Phi_{2}\;dx.
\end{eqnarray*}
To prove existence of weak solutions of the system $(\mathcal{S})$, we are going to study properties of the operators $L_1, L_2, B$ and $F$.\\
\textbf{1.} In view of the previous section in particular Lemmas \ref{thm6}, \ref{thm5},  \ref{thm7} and  we have similar properties to the operators $ L_1 $ and $ L_2$, i.e. $ L_1 $ and $ L_2$ are  demi-continuous and bounded, so their sum.\\
\textbf{2.} The second remark consist in the proof of coercivity of the operator $\widetilde{L}$ defined on the space $W_{p(x),q(x)}$ by: $(\widetilde{L}(u,v), (\Phi_{1},\Phi_{2}))= ((L_1-L_2+B)(u,v), (\Phi_{1},\Phi_{2})),$ for all $(\Phi_{1},\Phi_{2})\in W_{p(x),q(x)}$. Let $(u,v)\in W_{p(x),q(x)}$, then
\begin{eqnarray*}
(\widetilde{L }(u,v)  ,(u,v)) &\geq&\int_\Omega|\nabla
 u|^{p(x)}- \int_\Omega a(x)|u|^{p(x)} + \int_\Omega|\nabla  v|^{q(x)} -\int_\Omega c(x)|v|^{q(x)}\\
& & \quad +\; \int_\Omega  b(x)|u|^{\alpha(x)}|v|^{\beta(x) + 1} + \int_\Omega
 d(x)|u|^{\alpha(x) + 1}|v|^{\beta(x)}.
 \end{eqnarray*}
Since, the functionals $b(x)$ and $d(x)$ are positive on $\Omega$, we have
$$(\widetilde{L }(u,v)  ,(u,v)) \geq\int_\Omega|\nabla
 u|^{p(x)}- \int_\Omega a(x)|u|^{p(x)} + \int_\Omega|\nabla  v|^{q(x)} -\int_\Omega c(x)|v|^{q(x)}.$$
In view of inequality (\ref{e70}), we obtain
$$(\widetilde{L }(u,v)  ,(u,v))\geq  \Big(1-\frac{p^+}{p^-\lambda_{p }(a)}\Big)\int_\Omega|\nabla u|^{p(x)}
+ \Big(1-\frac{q^+}{q^-\lambda_{q }(c)}\Big)\int_\Omega|\nabla v|^{q(x)},$$
where $\lambda_{p }(a)$ and $\lambda_{q }(c)$ are respectively the first eigenvalue of the problem
$$\Delta_{p(x)}u=\lambda a(x) |u|^{p(x)-2}u\;\;\mbox{and}\;\;\Delta_{q(x)}u=\lambda c(x) |u|^{q(x)-2}u.$$
If we consider the fact that $\lambda_{p }(a)>\frac{p^+}{p^-}$ and $\lambda_{q }(c)>\frac{q^+}{q^-}$, we get
$$(\widetilde{L }(u,v)  ,(u,v))\geq  \int_\Omega|\nabla u|^{p(x)}
+ \int_\Omega|\nabla v|^{q(x)}.$$
Using inequalities (\ref{e3}) and (\ref{e4}), we obtain
$$(\widetilde{L }(u,v)  ,(u,v))\geq \min(|\nabla u|^{p^+}_{p(x)};|\nabla u|^{p^-}_{p(x)})+ \min(|\nabla v|^{q^+}_{q(x)};|\nabla v|^{q^-}_{q(x)}).$$
Since $\|u\|_{p(x)}=|\nabla u|_{p(x)}$, $\|v\|_{q(x)}=|\nabla v|_{q(x)}$ and $p^-, q^->1$, therefore
$$\frac{(\widetilde{L }(u,v)  ,(u,v))}{\|(u,v)\|_{p(x),q(x)}}\rightarrow\infty\;\;\mbox{ as}\;\;\|(u,v)\|_{p(x),q(x)}\rightarrow\infty.$$
The proof of the coercivity of the operator $\widetilde{L}$ is fulfilled.\\
\textbf{3.} The operator $B(u;v)$ is well defined; indeed, denotes
$$\Omega_1=\{x\in \Omega;\,|u(x)|\geq 1,\;|v(x)|\geq 1\},\hspace{0.5cm}\;\;\;\;\;\Omega_2=\{x\in \Omega;\,|u(x)|<1,\;|v(x)|<1\},$$
$$\Omega_3=\{x\in \Omega;\,|u(x)|\geq 1,\;|v(x)|\leq 1\}\;\;\mbox{and}\;\;\;\Omega_4=\{x\in \Omega;\,|u(x)|<1,\;|v(x)|\geq 1\}.$$
Clearly, we have
$$
\int_{\Omega} b(x)|u|^{\alpha(x)}|v|^{\beta(x)}v\phi_1 dx
=\sum_{i=1}^4\Big{(}\int_{\Omega_i} b(x)|u|^{\alpha(x)}|v|^{\beta(x)}vdx\phi_1\Big{)}.$$
Furthermore,
$$\Big|\int_{\Omega_1} b(x)|u|^{\alpha(x)}|v|^{\beta(x)}v\phi_1 dx\Big|\leq  \int_{\Omega_1} b(x)|u|^{\alpha^+}|v|^{\beta^++1}|\phi_1| dx.$$
Since $\alpha^++1<p^*(x),\;\beta^++1<q^*(x)$, then the following embeddings hold true
$$W_0^{1,p(x)}(\Omega)\hookrightarrow L^{\alpha^+p(x)}(\Omega)\;\mbox{and}\;\;W_0^{1,q(x)}(\Omega)\hookrightarrow L^{(\beta^++1)q(x)}(\Omega).$$
Then, we obtain
$$||u|^{\alpha^+}|_{\alpha^+p(x)}\leq c_1|u|_{p(x)}\leq c_2||u|^{\alpha^+}|_{p^*(x)},\;\;\mbox{and}\;\;||v|^{\beta^++1}|_{q(x)}\leq c_3||v|^{\beta^++1}|_{q^*(x)}.$$
If we apply (\ref{e3}), (\ref{e4})  and Proposition \ref{p22} and take the functionals $b(x)\in L^{s(x)}(\Omega);$  $d(x)\in L^{r(x)}(\Omega)$, then we have
\begin{equation*}
\Big|\int_{\Omega_1} b(x)|u|^{\alpha(x)}|v|^{\beta(x)}v\phi_1 dx\Big|\leq | b(x)|_{s(x)}|u^{\alpha^+ }|_{p^*(x)}||v|^{\beta^++1}|_{q^*(x)}|\phi_1|_{\widetilde{p}(x)}<\infty.
\end{equation*}
\begin{equation*}
\Big|\int_{\Omega_1} d(x)|v|^{\beta(x)}|u|^{\alpha(x)}u\phi_2 dx\Big|\leq | d(x)|_{r(x)}||u|^{\alpha^++1}|_{p^*(x)}
||v|^{\beta^+}|_{q^*(x)}|\phi_2|_{\widetilde{q}(x)}<\infty.
\end{equation*}
\begin{equation*}
\Big|\int_{\Omega_2} b(x)|u|^{\alpha(x)}|v|^{\beta(x)}v\phi_1 dx\Big|\leq | b(x)|_{s(x)}|u^{\alpha^- }|_{p^*(x)}||v|^{\beta^-+1}|_{q^*(x)}|\phi_1|_{\widetilde{p}(x)}<\infty.
\end{equation*}
\begin{equation*}
\Big|\int_{\Omega_2} d(x)|v|^{\beta(x)}|u|^{\alpha(x)}u\phi_2 dx\Big|\leq | d(x)|_{r(x)}||u|^{\alpha^-+1}|_{p^*(x)}
||v|^{\beta^-}|_{q^*(x)}|\phi_2|_{\widetilde{q}(x)}<\infty.
\end{equation*}
Repeating the same arguments we deduce
\begin{equation*}
\Big|\int_{\Omega_i} b(x)|u|^{\alpha(x)}|v|^{\beta(x)}v\phi_1 dx\Big|<\infty,\;
\Big|\int_{\Omega_i} d(x)|v|^{\beta(x)}|u|^{\alpha(x)}u\phi_2 dx\Big|<\infty,\;\mbox{for}\;i=3,4.
\end{equation*}
Hence, $|(B(u;v),(\Phi_{1},\Phi_{2}))|<\infty.$  The operator $B(u;v)$ is well defined on $W_{p(x),q(x)}$.
\end{Remark}
\begin{Proof}
Using remark \ref{R2} and Proposition \ref{p0} it remaind to prove the continuity of the operator $B$. To this end we will show the compactness of $B$.\\
Let $\{(u_n, v_n)\}\subset W_{p(x),q(x)}$ be a sequence such that $(u_n, v_n)\rightharpoonup (u, v)$ weakly in $W_{p(x),q(x)}$. We claim that $B(u_n, v_n)\rightarrow  B(u, v)$ strongly in $W_{p(x),q(x)}$, i.e. for all $(\Phi_{1},\Phi_{2})\in W_{p(x),q(x)}$ we have
$$\Big{|}(B(u_n, v_n)- B(u, v);\;(\Phi_{1},\Phi_{2}))\Big{|}=\circ(1)\;\;\mbox{as}\;\;n\,\rightarrow\infty.$$
Clearly
$$B(u_n, v_n)- B(u, v)=(B_u(u_n, v_n)- B_u(u, v))+(B_v(u_n, v_n)- B_v(u, v)),$$
where
$$(B_u(u_n, v_n)- B_u(u, v);\;(\Phi_{1},\Phi_{2}))=\int_\Omega b(x)(|u_n|^{\alpha(x)}|v_n|^{\beta(x)}v_n-|u|^{\alpha(x)}|v|^{\beta(x)}v)\Phi_{1}dx,$$
and
$$(B_v(u_n, v_n)- B_v(u, v);\;(\Phi_{1},\Phi_{2}))=\int_\Omega d(x)(|v_n|^{\beta(x)}|u_n|^{\alpha(x)}u_n-|v|^{\beta(x)}|u|^{\alpha(x)}u)\Phi_{2}dx.$$
Then it's sufficient to prove the compactness of  $B_u(u,v)$ and $B_v(u,v).$
\begin{eqnarray*}
(B_u(u_n, v_n)- B_u(u, v);\;(\Phi_{1},\Phi_{2})) & = & \int_\Omega b(x)|v_n|^{\beta(x)+1}(|u_n|^{\alpha(x)}-|u|^{\alpha(x)})\Phi_{1}\,dx\\
& & +\int_\Omega b(x)|u|^{\alpha(x)}(|v_n|^{\beta(x)+1}-|v|^{\beta(x)}v)\Phi_{1}\,dx.
\end{eqnarray*}
In view of Remark \ref{R2}, precisely item \textbf{3.} one writes
\begin{eqnarray*}
\Big{|}(B_u(u_n, v_n)- B_u(u, v);\;(\Phi_{1},\Phi_{2}))\Big{|} & \leq & c_1 |b(x)|_{s(x)}\Big{(}|v_n|^{\beta(x)+1}_{q^*(x)}\Big{|}|u_n|^{\alpha(x)}-
|u|^{\alpha(x)}\Big{|}_{p^*(x)}\\
& &||u|^{\alpha(x)}|_{p^*(x)}|
|v_n|^{\beta(x)+1}-|v|^{\beta(x)}v|_{q^*(x)}\Big{)}|\Phi_{1}|_{\widetilde{p}(x)}.
\end{eqnarray*}
Similar calculation gives us the following inequality
\begin{eqnarray*}
\Big{|}(B_v(u_n, v_n)- B_v(u, v);\;(\Phi_{1},\Phi_{2}))\Big{|} & \leq & c_2 |d(x)|_{r(x)}\Big{(}|u_n|^{\alpha+1}_{p^*(x)}\Big{|}|v_n|^{\beta(x)}-
|v|^{\beta(x)}\Big{|}_{q^*(x)}\\
& &||v|^{\beta}|_{p^*(x)}|
|u_n|^{\alpha(x)+1}-|u|^{\alpha(x)}u|_{p^*(x)}\Big{)}|\Phi_{2}|_{\widetilde{q}(x)}.
\end{eqnarray*}
Due to the continuity of Nemytskii operators $u\rightarrow |u|^{\alpha(x)}$ $($resp. $v\rightarrow |v|^{\beta(x)} v)$ from $L^{p(x)}(\Omega)$ into $L^{p^*(x)}(\Omega)$ $($resp.
from $L^{q(x)}(\Omega)$ into $L^{q^*(x)}(\Omega)),$ there exists $n_0\geq 0$ such that for all $n\geq n_0$ we have
\begin{equation}\label{e142}
\Big{|}|u_n|^{\alpha(x)}-
|u|^{\alpha(x)}\Big{|}_{p^*(x)}=\circ(1),
\end{equation}
\begin{equation}\label{e143}
\Big{|}|v_n|^{\beta(x)+1}-|v|^{\beta(x)}v\Big{|}_{q^*(x)}=\circ(1).
\end{equation}
Finally from equations (\ref{e142}) and (\ref{e143}), we have the claim and the operator $B$ will be compact and completely continuous. Hence, $B$ satisfies the $M_0$-condition and the system $(\mathcal{S})$ possess a weak solution $(u,v)\in W_{p(x),q(x)}$, for all $(f,g)$ in the dual of  $W_{p(x),q(x)}$. The proof of the main result on bounded domains is completed.
\end{Proof}
\section{Nonlinear systems defined on $\mathbb{R}^N$}
In this section, we study existence of weak solutions of the following system.
\begin{eqnarray} \label{S}
\begin{gathered}
-\Delta_{p(x)}u = a(x)|u|^{p(x)-2}u-b(x)|u|^{\alpha(x)}|v|^{\beta(x)} v+f,\\
-\Delta_{q(x)}v = c(x)|u|^{\alpha(x)} |v|^{\beta(x)} u - d(x) |v|^{q(x)-2}v +g ,\\
\lim_{|x|\to\infty}u=\lim_{|x|\to\infty}v=0\quad u,v>0
\end{gathered}
\end{eqnarray}
which is defined on $\mathbb{R}^N$.
We assume that the coefficients
$a(x), b(x), c(x), d(x)$ are smooth positive functions satisfying assumptions F1) and F2) introduced in section 2. In addition, functionals $\alpha(x)$ and $\beta(x)$ will be such that condition F3). Note that we conserve notations of section 4 with $\Omega=\mathbb{R}^N$, in particular $W_{p(x),q(x)}$ represent the product space $W_0^{1,p(x)}(\mathbb{R}^N)\times W_0^{1,q(x)}(\mathbb{R}^N)$. By transforming the weak formulation for the system (\ref{S}) to the operator formulation, we will get the same operators $\widetilde{L}, L_1, L_2, B$ and $F$ which take similar definitions in Remark \ref{R2}.
\begin{Remark}\label{R3}
It's well known that the operator $L_1$ is well defined, continuous on $W_{p(x),q(x)}$, for the proof we cite the work of \cite{O}.
\end{Remark}
\begin{Lemma}\label{l5}
Under the assumptions F1), F2) and F3), The operators $L_2$ and $ B$ are well defined on the space $W_{p(x),q(x)}$.
\end{Lemma}
\begin{Proof}
For all pairs of real functions $(u,v), (\Phi_1, \Phi_2)\in W_{p(x),q(x)}$, under the assumptions F1), F2) and F3), we can write
\begin{eqnarray*}
\Big{|}(L_2 (u,v),(\Phi_1, \Phi_2))\Big{|}&=&\Big{|}\int_{\mathbb{R}^N} a(x)
|u|^{p(x)-2}u\Phi_{1}dx+\int_{\mathbb{R}^N} c(x) |v|^{q(x)-2}v \Phi_{2}dx\Big{|}\\
& \leq & \int_{\mathbb{R}^N} a(x)
|u|^{p(x)-1}|\Phi_{1}|dx+\int_{\mathbb{R}^N} c(x) |v|^{q(x)-1} |\Phi_{2}|dx\\
& \leq & \int_{\mathbb{R}^N} a(x)
|u|^{p^+-1}|\Phi_{1}|dx+\int_{\mathbb{R}^N} c(x) |v|^{q^+-1} |\Phi_{2}|dx.
\end{eqnarray*}
If we consider the fact that
$$W^{1,p(x)}_{0}(\Omega)\hookrightarrow L^{(p^+-1)p(x)}(\Omega) \Rightarrow||u|^{p^+-1}|_{p(x)}=|u|^{p^+-1}_{(p^+-1)p(x)}\leq c\|u|^{p^+-1}|_{(p^+-1)p(x)},$$
and if we apply (\ref{e3}), (\ref{e4}), Proposition \ref{p22} and take $a(x)\in L^{k_1(x)}(\mathbb{R}^N)$, $c(x)\in L^{k_2(x)}(\mathbb{R}^N)$ then we have
\begin{eqnarray*}
\Big{|}(L_2 (u,v),(\Phi_1, \Phi_2))\Big{|}&\leq & c\Big{(}|a(x)|_{k_1(x)}||u|^{p^+-1}|_{p^*(x)}|\Phi_{1}|_{\widetilde{p}(x)}\\
&  &\hspace{2cm}+|c(x)|_{k_2(x)}||v|^{q^+-1}|_{q^*(x)}
|\Phi_{2}|_{\widetilde{q}(x)}\Big{)}.
\end{eqnarray*}
Therefore, the operator $L_2$ is well defined. Note that $B$ is again well defined on  $W_{p(x),q(x)}$, the proof is the same as in the Remark \ref{R2}, item \textbf{3.} replacing $\Omega$ by $\mathbb{R}^N$.
\end{Proof}
Next, we deal with the demicontinuity of the operator $B$. For this aim, we denote by $B_r$ the ball of radius $r$ which is centered at the origin of $\mathbb{R}^N$ and let $B'_r$ the complementary of  $B_r$ in $\mathbb{R}^N$. i.e. $B'_r=\mathbb{R}^N-B_r$.
\begin{Lemma}\label{l6}
Under the assumptions F1), F2) and F3), The operators $B$ is demi-continuous on the space $W_{p(x),q(x)}$.
\end{Lemma}\label{l21}
\begin{Proof} Recall that for all pairs of real functions $(u,v), (\Phi_1, \Phi_2)\in W_{p(x),q(x)}$, 
$$(B(u,v), (\Phi_{1},\Phi_{2}))=
 \int_{\mathbb{R}^N} b(x) |u|^{\alpha(x)} |v|^{\beta(x)} v \Phi_{1}dx+
\int_{\mathbb{R}^N} d (x) |v|^{\beta(x)} |u|^{\alpha(x)} u \Phi_{2}dx.$$
Let $\{(u_n, v_n)\}\subset W_{p(x),q(x)}$ be a sequence such that $(u_n, v_n)\rightharpoonup (u, v)$ weakly in $W_{p(x),q(x)}$. We claim that $B(u_n, v_n)\rightarrow  B(u, v)$ strongly in $W_{p(x),q(x)}$.\\ 
Repeating calculations of the proof of Theorem 2, one writes\\
$$B(u_n, v_n)- B(u, v)=(B_u(u_n, v_n)- B_u(u, v))+(B_v(u_n, v_n)- B_v(u, v)),$$
where
$$(B_u(u_n, v_n)- B_u(u, v);\;(\Phi_{1},\Phi_{2}))=\int_{\mathbb{R}^N} b(x)(|u_n|^{\alpha(x)}|v_n|^{\beta(x)}v_n-|u|^{\alpha(x)}|v|^{\beta(x)}v)\Phi_{1}dx,$$
$$(B_v(u_n, v_n)- B_v(u, v);\;(\Phi_{1},\Phi_{2}))=\int_{\mathbb{R}^N} d(x)(|v_n|^{\beta(x)}|u_n|^{\alpha(x)}u_n-|v|^{\beta(x)}|u|^{\alpha(x)}u)\Phi_{2}dx.$$
On the other hand, we have
\begin{eqnarray*}
(B_u(u_n, v_n)- B_u(u, v);\;(\Phi_{1},\Phi_{2}))&=& \int_{B_r} b(x)(|u_n|^{\alpha(x)}|v_n|^{\beta(x)}v_n-|u|^{\alpha(x)}|v|^{\beta(x)}v)\Phi_{1}\\
&+&\int_{B'_r} b(x)(|u_n|^{\alpha(x)}|v_n|^{\beta(x)}v_n-|u|^{\alpha(x)}|v|^{\beta(x)}v)\Phi_{1}.\\
&=&I_1+I_2.
\end{eqnarray*}
Since $B_r$ is a bounded domain in $\mathbb{R}^N$. Using the proof of Theorem 2, we obtain the demicontinuity of the operator $B$ on the unit ball $B_r$. Hence it remainds to justify the demicontinuity of $B$ on $B'_r$.\\
In view of Remark \ref{R2}, precisely item \textbf{3.} one writes
\begin{eqnarray*}
\Big{|}I_2\Big{|}&\leq & c |b(x)|_{s(x)}\Big{(}|v_n|^{\beta(x)+1}_{q^*(x)}\Big{|}|u_n|^{\alpha(x)}\Big{)}\\
&&-c |b(x)|_{s(x)}\Big{(}
|u|^{\alpha(x)}\Big{|}_{p^*(x)}
||u|^{\alpha(x)}|_{p^*(x)}|
|v_n|^{\beta(x)+1}-|v|^{\beta(x)}v|_{q^*(x)}\Big{)}|\Phi_{1}|_{\widetilde{p}(x)}.
\end{eqnarray*} 
Due to the continuity of Nemytskii operators $u\rightarrow |u|^{\alpha(x)}$ $($resp. $v\rightarrow |v|^{\beta(x)} v)$  and the fact that
$|b(x)|_{L^{s(x)}(B'_r)}\rightarrow 0\;\;\mbox{for}\;\;r\rightarrow \infty,$
It follows that the operator $B_u$ satisfies the $M_0$-condition.
Similarly, the operator $B_v$ will be demicontinuous and  so the operator $B$. That's completes the proof of Lemma \ref{l21}.\\
\end{Proof}
\begin{Remark}
\textbf{1.} The proof of coercivity of the operator $\widetilde{L}$ is similar to each in bounded domains.\\
\textbf{2.} All conditions of Proposition \ref{p0} are satisfied by the operator $\widetilde{L}$ on $\mathbb{R}^N$, which guarantees the  existence of a weak solution for system (\ref{S}).
\end{Remark}

\end{document}